\documentstyle[12pt,amsfonts,amssymb,leqno]{article}

\newtheorem{proposition}{Proposition}[section]
\newtheorem{thm}[proposition]{Theorem}

\newtheorem{lemma}[proposition]{Lemma}
\newtheorem{defn}[proposition]{Definition}

\newcommand{\zerohandgrenade}{0^{\: \mid^{\! \! \! \bullet}}} 

\begin{document}

\begin{center}
{\Large \bf Sharps and the
$\Sigma^1_3$ correctness of $K$} 
\end{center}
\begin{center}
\renewcommand{\thefootnote}{\fnsymbol{footnote}}
\renewcommand{\thefootnote}{arabic{footnote}}
\renewcommand{\thefootnote}{\fnsymbol{footnote}}
{\large Ralf Schindler}
\renewcommand{\thefootnote}{arabic{footnote}}
\end{center}
\begin{center} 
{\footnotesize
{\it Institut f\"ur Formale Logik, Universit\"at Wien, 1090 Wien, Austria}} 
\end{center}

\begin{center}
{\tt rds@logic.univie.ac.at}

{\tt http://www.logic.univie.ac.at/${}^\sim$rds/}\\
\end{center}

The purpose of the present paper is to present a new, simple, and purely combinatorial 
proof of the following result.

\begin{thm}\label{thm} {\em \bf (Steel-Welch 1993, \cite[Theorem 4.1]{johnphilip} )} 
Let $A \subset {\mathbb R}$ be $\Pi^1_2$.
Suppose that there is some sequence $(x_n \colon n<\omega)$ such that
$x_0 \in A$ and for all $n<\omega$, $x_{n+1} = x_n^\#$. Suppose also
that there is some $N < \omega$ such that there is no inner model with $N$ strong
cardinals.
Then $A \cap K \not= \emptyset$.
\end{thm}

Here, $K$ denotes the core model; cf.~the remark right after Definition \ref{K}.
It is open whether Theorem \ref{thm} still holds if we replace the second sentence in
its statement by ``Suppose that there is some $x \in A$ such that $x^\#$ exists.'' It
is also open 
whether Theorem \ref{thm} still holds if we replace the third sentence in
its statement by ``Suppose also that there is no inner model with a Woodin cardinal,
but $K$ exists'' (cf.~Definition \ref{K}).

We refer the reader to \cite{correctness}. The current argument will exploit, among
other things,
the argument of \cite{correctness}. 

\begin{defn}\label{K}
Let ${\cal A}$ be a transitive model of ${\sf ZFC}$. Then by $K^{\cal A}$ we denote
the model which is recursively constructed inside ${\cal A}$ in the manner of
\cite[\S 6]{CMIP}, if it exists (otherwise we let $K^{\cal A}$ undefined). 
If $K^{\cal A} \downarrow$ then we say that $K^{\cal A}$ exists.
If $K^V \downarrow$ then we write $K = K^V$ and say that $K$ exists. 
\end{defn}

It is shown in \cite{habil} that $K$ exists if
$\zerohandgrenade$ doesn't exist.
It is also shown in \cite{CMIP} that $K^M$ exists if $M = V_\Omega^{\cal H}$, where
${\cal H}$ is a transitive model of ``${\sf ZFC}^- + \Omega$ is measurable + there is
no inner model with a Woodin cardinal'' (in this case
we'll sometimes also write $K^{\cal H}$ for
$K^M$).

We shall prove Theorem \ref{thm} with the third sentence in its statement 
being replaced by ``Suppose also that $0^\P$ doesn't exist.'' We'll leave a proof of
Theorem \ref{thm} as stated as an exercise to the reader.

\begin{defn}\label{P}
Let $(\clubsuit)$ denote the following assertion.  
Let $x \in {\mathbb R}$ be such that $x^\#$
exists. If $K^{L[x]}$ and $K^{L[x^\#]}$ both exist and are coiterable
then there is some $\alpha \in
{\rm OR}$ such that $K^{L[x^\#]}||\alpha$ iterates past $K^{L[x]}$.
\end{defn}

\begin{lemma}\label{lemma1}
Suppose 
that $0^\P$ doesn't exist. Then $(\clubsuit)$ holds.
\end{lemma}

{\sc Proof.} Suppose that $x \in {\mathbb R}$ witnesses the failure of $(\clubsuit)$.
It is fairly easy too see that then $K^{L[x^\#]} \models$ ``there is a strong
cardinal.'' (Cf.~the proof of \cite[Lemma 3.3]{johnphilip}.) Let us now assume that 
$0^\P$ doesn't exist. We aim to derive a contradiction. Let us work in $L[x^\#]$.

\bigskip
{\bf Claim 1.} ${\rm cf}(\kappa^{+K^{L[x]}}) = \omega$ for all $\kappa$.

\bigskip
{\sc Proof.} It is true that ${\rm cf}(\kappa^{+{L[x]}}) = \omega$ for all $\kappa$.
Let us thus fix some $K^{L[x]}$-cardinal 
$\kappa$ such that $\kappa^{+K^{L[x]}} < \kappa^{+{L[x]}}$. 
Let $\lambda = {\rm Card}^{L[x]}(\kappa) \leq \kappa$.
Then
$\lambda$ is neither an $x$-indiscernible nor singular in $L[x]$. Let $\eta < \lambda$
be the largest $x$-indiscernible which is smaller than $\lambda$.

Let $\tau_n$ enumerate the Skolem terms of $L[x]$. The sequence $(\lambda_n \colon
n<\omega)$, where 
$$\lambda_n = {\rm sup}( \{ \tau_n({\vec \xi}) \colon {\vec \xi} < \eta \} \cap
\lambda )
< \lambda$$ 
witnesses that ${\rm cf}(\lambda) = \omega$. But we'll have that
${\rm cf}^{L[x]}(\kappa^{+K^{L[x]}}) =
\lambda$ by weak covering applied inside $L[x]$ (cf.~\cite{covering}).

\hfill $\square$ (Claim 1)

\bigskip 
{\bf Claim 2.} $K^{L[x]}$ doesn't move in the comparison with $K^{L[x^\#]}$.

\bigskip
{\sc Proof sketch.} The point is that by our assumption $K^{L[x]}$ absorbs all
coiterable set-sized premice which exist in $L[x^\#]$. Jensen's argument yielding that 
below $0^\P$ any
universal weasel is an iterate of $K$ then gives this Claim.

\hfill $\square$ (Claim 2)

\bigskip
It is now easy to see that Claims 1 and 2, combined with an application of 
weak covering applied inside $L[x^\#]$ (cf.~\cite{covering}) yields the Lemma.

\hfill $\square$ (Lemma \ref{lemma1})
 
\bigskip
The proof of Lemma \ref{lemma1} is certainly more interesting than its result. If we
had assumed the existence of $x^{\#\#}$ then we could have just cited 
\cite[Lemma 3.3]{johnphilip}. We conjecture that
$(\clubsuit)$ still holds under much weaker assuptions 
than the non-existence of $0^\P$ (cf.~\cite[p.~188, Question3]{johnphilip}). 

We are now going to prove the following result, which will
immediately give Theorem \ref{thm} (the third sentence in its statement 
being replaced by ``Suppose also that $0^\P$ doesn't exist'') via Lemma \ref{lemma1}.

We emphasize that Theorem \ref{thm1} is not given by the results of \cite{johnphilip};
the proof of \cite[Theorem 4.1]{johnphilip} which is given in \cite{johnphilip}
heavily uses universal iterations which
are not known to exist significantly above $0^\P$.

\begin{thm}\label{thm1} 
Let $A \subset {\mathbb R}$ be $\Pi^1_2$.
Suppose that there is some sequence ${\vec x} = (x_n \colon n<\omega)$ such that
$x_0 \in A$ and for all $n<\omega$, $x_{n+1} = x_n^\#$. Suppose also
that there is no inner model with a Woodin cardinal, that $K^{L[{\vec x}]}$ exists, and
that $(\clubsuit)$ holds.
Then $A \cap K^{L[{\vec x}]} 
\not= \emptyset$.
\end{thm}

\bigskip
{\sc Proof} of Theorem \ref{thm1}. Let $A = \{ z \in {\mathbb R} \colon \Phi(z)
\}$ where $\Phi(-)$ is $\Pi^1_2$. There is a tree $T \in K^{L[{\vec x}]}$ searching
for a quadruple $({\vec y},{\vec M},{\vec {\cal T}},{\vec \sigma})$ such that the
following hold true.

\bigskip
\noindent $\bullet \ \ $ ${\vec y} = (y_n \colon n<\omega) \in {}^\omega {\mathbb R}$, 

\noindent $\bullet \ \ $ ${\vec M} = (M_n \colon n<\omega)$ such that 
for all $n<\omega$ do we have the following: 

(a) $M_n = (J_{\alpha_n}[y_n];\in,y_n,U_n)$ for some $\alpha_n$, $U_n$, 

(b) 
$M_{n+1} \models$ ``$M_n = y_n^\#$'' (in particular, $M_{n+1}$ thinks that $M_n$ is
iterable), 

(c) $y_{n+1}$ is the master code of $M_n$, 

(d) $M_0 \models \Phi(y_0)$,
and 

(e) setting $\kappa = {\rm crit}(U_{n+1})$, there is an initial segment of
$K^{M_{n+1}}$ which iterates past $K^{L_{\kappa}[y_n]}$, and

\noindent $\bullet \ \ $ $({\vec {\cal T}},{\vec \sigma})$ witnesses that each 
individual 
$K^{M_{n}}$ ($n<\omega)$ is iterable (cf.~\cite{correctness}), i.e., ${\vec T} =
({\cal T}_n
\colon n<\omega)$, ${\vec \sigma} = (\sigma_n \colon n<\omega)$, and
for all $n<\omega$ do we have the following:

(a) ${\cal T}_n$ is a countable tree of successor length on $K^{L[{\vec x}]}$, and

(b) $\sigma_n \colon K^{M_n} \rightarrow {\cal M}^{{\cal T}_n}_\infty || \beta_n$,
some $\beta_n \leq {\cal M}^{{\cal T}_n}_\infty \cap {\rm OR}$, is elementary. 

\bigskip
We are now going to prove that $$\emptyset \not= p[T] = \{ y_0 \colon \exists
(y_1,y_2,...) \exists {\vec M} \exists {\vec {\cal T}} \exists {\vec \sigma}
((y_0, y_1, ...),{\vec M},{\vec {\cal T}},{\vec \sigma}) \in [T] \} \subset A.$$
We may well leave the verification of $p[T] \not= \emptyset$ as an exercise to the
reader. 

Now fix $({\vec y},{\vec M},{\vec {\cal T}},{\vec \sigma}) \in T$. Let ${\vec y} =
(y_n \colon n<\omega)$ and ${\vec M} = (M_n \colon n<\omega)$. Let us prove that
$y_0 \in A$. Let $({}^\alpha_n)$ denote the assertion that the $\alpha^{\rm th}$
iterate of $M_n$ is well-founded. It clearly suffices to prove the following.

\bigskip
{\bf Main Claim.} For all $\alpha$, for all $n$, $({}^\alpha_n)$ holds.

\begin{defn}\label{ui}
Let $n<\omega$. We write $(M_n^i, \pi_n^{ij} \colon i \leq j \leq
\alpha)$ for the putative iteration of $M_n$ of length $\alpha +1$, if it exists; and
if so then
for $i<\alpha$ we write $\kappa_n^i$ for the critical
point of $\pi_n^{0i}(U_n)$, i.e., of the top extender of $M_n^i$.
We call $\alpha$ a {\em uniform indiscernible} provided that
for all $n<\omega$, the putative iteration $(M_n^i, \pi_n^{ij} \colon i \leq j \leq
\alpha)$ of $M_n$ of length $\alpha +1$ exists and 
$\{ \kappa_i \colon i<\alpha \}$ is (closed and) unbounded in
$\alpha$.
\end{defn}

{\sc Proof} of the Main Claim. We'll prove the Main Claim by induction on $\alpha$. 

\bigskip
{\sc Case 1.} $\alpha$ is not a uniform indiscernible.

\bigskip
Let $n<\omega$. 
%Let $(M_{n+1}^i, \pi_{n+1}^{ij} \colon i \leq j \leq
%\alpha)$ be the putative iteration of $M_{n+1}$ of length $\alpha +1$, which exists by
%our inductive hypothesis. 
By our case assumption, there are some $m > n$ and 
$\beta < \alpha$ such
that $\alpha \in M_m^\beta$. But $M_m^\beta \models$ ``$M_n = y_n^\#$,'' so
that we may argue inside $M_m^\beta$ and deduce that the $\alpha^{\rm th}$ iterate
of $M_n$, viz.~$M_n^\alpha$, is well-founded.

\bigskip
{\sc Case 2.} $\alpha$ is a uniform indiscernible.

\bigskip
Let $n<\omega$.
Let $\kappa = {\rm crit}(U_{n+1}) = \kappa_{n+1}^0$, and let 
${\cal P}$ be the proper 
initial segment of $K^{M_{n+1}} = K^{L_\kappa[y_{n+1}]}$ 
which iterates past $K^{L_\kappa[y_n]}$.
Let
$({\cal T},{\cal U})$ be the coiteration of 
$K^{L_\kappa[y_n]}$ with $K^{L_\kappa[y_{n+1}]}$. 
$(M_n^i, \pi_n^{ij} \colon i \leq j <
\alpha)$ is the putative iteration of $M_n$ of length $\alpha +1$.
Let $$\sigma \colon M_{n+1} \rightarrow_{U_{n+1}} M' {\rm , }$$
i.e., $\sigma = \pi_{n+1}^{01}$ and $M' = M_{n+1}^1$.

Let $X \in {\cal P}(\kappa) \cap {\cal M}_\infty^{\cal U}$. Then $X
= \pi_{i \infty}^{\cal U}({\bar X})$, some $i < \kappa$, ${\bar X}$, and 
$\sigma(X) = \pi_{i \infty}^{\sigma({\cal U})}({\bar X}) = \pi_{\kappa
\infty}^{\sigma({\cal U})}(X)$. Therefore, $\pi_{\kappa
\infty}^{\sigma({\cal U})} \upharpoonright \kappa^{+{\cal M}_\infty^{\cal U}} =
\sigma \upharpoonright \kappa^{+{\cal M}_\infty^{\cal U}}$. The same argument shows
that $\pi_n^{\kappa \sigma(\kappa)} \upharpoonright \kappa^{+{M_n^\kappa}} =
\sigma \upharpoonright \kappa^{+{M_n^\kappa}}$; we construe $M_n$ in such a way that
${\rm crit}(U_n)^{+M_n} = M_n \cap {\rm OR}$, so that this latter equality means that
$\pi_n^{\kappa \sigma(\kappa)} \upharpoonright M_n^\kappa \cap {\rm OR} =
\sigma \upharpoonright M_n^\kappa \cap {\rm OR}$.
Let us write $\lambda = \sigma(\kappa) = \kappa_{n+1}^1$.

Now $(\sigma({\cal T}),\sigma({\cal U}))$ is the coiteration of 
$K^{L_{\lambda}[y_n]}$ with $K^{L_{\lambda}[y_{n+1}]}$.
We'll have that $\kappa^{+{\cal M}_\infty^{\cal U}} = 
\kappa^{+{\cal M}_\kappa^{\sigma({\cal U})}} = 
\kappa^{+{\cal M}_\kappa^{\sigma({\cal T})}} \geq
\kappa^{+K^{L_{\lambda}[y_n]}} = \kappa^{+{L_{\lambda}[y_n]}}
= \kappa^{+M_n^\kappa}$, so that we get that
$$\pi_n^{\kappa \lambda} \upharpoonright M_n^\kappa \cap {\rm OR} = \pi_{\kappa
\infty}^{\sigma({\cal U})} \upharpoonright M_n^\kappa \cap {\rm OR}.$$

This now buys us that if we let $({\cal T}^*,{\cal U}^*)$ denote the coiteration of 
$K^{L_\alpha[y_n]}$ with $K^{L_\alpha[y_{n+1}]}$
then for typical $i \leq j < \alpha$ (namely, for all $i \leq j \in \{
\kappa_{n+1}^\beta \colon \beta < \alpha \}$) 
we'll have that
$$\pi_n^{ij} \upharpoonright M_n^i \cap {\rm OR} =
\pi_{i j}^{{\cal U}^*} 
\upharpoonright M_n^i \cap {\rm OR} {\rm . }$$
Moreover, ${\cal U}^*$ may be construed as an iteration of ${\cal P}$.
As ${\cal P}$ is iterable, we may thus conclude that the $\alpha^{\rm th}$ 
iterate 
of $M_n$, viz. $M_n^\alpha$, is well-founded (cf.~the argument of \cite{correctness}).

\bigskip
\hfill $\square$ (Main Claim)

\hfill $\square$ (Theorem \ref{thm1})

\bigskip
Using \cite{bill} it can be verified that Theorem \ref{thm1} still holds if the
assumption that $K^{L[{\vec x}]}$ is being crossed out and the conclusion is being
replaced by ``Then there is some lightface iterable premouse ${\cal M}$ with $A \cap
{\cal M} \not= \emptyset$.

\end{document}